\documentclass{amsart}%
\usepackage{amssymb}
\usepackage{amsmath}
\usepackage{amsfonts}%
\usepackage{graphicx}
\newtheorem{theorem}{Theorem}
\theoremstyle{plain}

\newtheorem{corollary}{Corollary}

\newtheorem{lemma}{Lemma}

\newtheorem{proposition}{Proposition}

\numberwithin{equation}{section}
\begin{document}
\title[Very well-covered graphs with log-concave independence polynomials]{Very well-covered graphs with log-concave independence polynomials}
\author{Vadim E. Levit and Eugen Mandrescu}
\address{Department of Computer Science\\
Holon Academic Institute of Technology\\
52 Golomb Str., Holon 58102, ISRAEL}
\email{\{levitv, eugen\_m\}@hait.ac.il}
\date{November 10, 2004}

\begin{abstract}
If $s_{k}$ equals the number of stable sets of cardinality $k$ in the graph
$G$, then $I(G;x)=\sum\limits_{k=0}^{\alpha(G)}s_{k}x^{k}$ is the
\textit{independence polynomial} of $G$ (Gutman and Harary, 1983). Alavi,
Malde, Schwenk and Erd\"{o}s (1987) conjectured that $I(G;x)$ is unimodal
whenever $G$ is a forest, while Brown, Dilcher and Nowakowski (2000)
conjectured that $I(G;x)$ is unimodal for any well-covered graph $G$. Michael
and Traves (2003) showed that the assertion is false for well-covered graphs
with $\alpha(G)\geq4$, while for very well-covered graphs the conjecture is
still open.

In this paper we give support to both conjectures by demonstrating that if
$\alpha(G)\leq3$, or $G\in\{K_{1,n},P_{n}:n\geq1\}$, then $I(G^{\ast};x)$ is
log-concave, and, hence, unimodal (where $G^{\ast}$ is the very well-covered
graph obtained from $G$ by appending a single pendant edge to each vertex).

\end{abstract}
\maketitle

\section{Introduction}

Throughout this paper $G=(V,E)$ is a finite, undirected, loopless and without
multiple edges graph with vertex set $V=V(G)$ and edge set $E=E(G)$. The set
$N(v)=\{u:u\in V,uv\in E\}$ is the \textit{neighborhood} of $v\in V$, and
$N[v]=N(v)\cup\{v\}$. As usual, a \textit{tree} is an acyclic connected graph,
while a \textit{spider} is a tree having at most one vertex of degree $\geq3$.
$K_{n},P_{n},K_{n_{_{1}},n_{_{2}},...,n_{p}}$ denote, respectively, the
complete graph on $n\geq1$ vertices, the chordless path on $n\geq1$ vertices,
and the complete $p$-partite graph on $n_{1}+n_{2}+...+n_{p}$ vertices,
$n_{1},n_{2},...,n_{p}\geq1$. A graph is called \textit{claw-free} if it has
no induced subgraph isomorphic to $K_{1,3}$. The \textit{disjoint union} of
the graphs $G_{1},G_{2}$ is the graph $G=G_{1}\sqcup G_{2}$ having
$V(G)=V(G_{1})\cup V(G_{2})$ and $E(G)=E(G_{1})\cup E(G_{2})$. If $G_{1}%
,G_{2}$ are disjoint graphs, then their \textit{Zykov sum}, (\cite{Zykov}), is
the graph $G_{1}\uplus G_{2}$ with $V(G_{1}\uplus G_{2})=V(G_{1})\cup
V(G_{2})$ and $E(G_{1}\uplus G_{2})=E(G_{1})\cup E(G_{2})\cup\{v_{1}%
v_{2}:v_{1}\in V(G_{1}),v_{2}\in V(G_{2})\}$. In particular, $\sqcup nG$ and
$\uplus nG$ denote the disjoint union and Zykov sum, respectively, of $n>1$
copies of the graph $G$.

A \textit{stable} set in $G$ is a set of pairwise non-adjacent vertices. The
\textit{stability number} $\alpha(G)$ of $G$ is the maximum size of a stable
set in $G$. A graph $G$ is called \textit{well-covered} if all its maximal
stable sets are of the same cardinality, \cite{Plum}. If, in addition, $G$ has
no isolated vertices and its order equals $2\alpha(G)$, then $G$ is
\textit{very well-covered}, \cite{Fav1}. By $G^{*}$ we mean the graph obtained
from $G$ by appending a single pendant edge to each vertex of $G$. Let us
remark that $G^{*}$ is well-covered (see, for instance, \cite{LevMan0}), and
$\alpha(G^{*})=n$. In fact, $G^{*}$ is very well-covered.

Let $s_{k}$ be the number of stable sets in $G$ of cardinality $k\in
\{0,1,...,\alpha(G)\}$. The polynomial $I(G;x)=\sum\limits_{k=0}^{\alpha
(G)}s_{k}x^{k}=1+s_{1}x+s_{2}x^{2}+...+s_{\alpha}x^{\alpha},\alpha=\alpha(G)$
is called the \textit{independence polynomial} of $G$ (Gutman and Harary,
\cite{GutHar}). In \cite{GutHar} was also proved the following equalities.

\begin{proposition}
\label{prop4}If $v\in V(G)$, then $I(G;x)=I(G-v;x)+xI(G-N[v];x)$, and
\[
I(G_{1}\sqcup G_{2};x)=I(G_{1};x)\cdot I(G_{2};x),\quad I(G_{1}\uplus
G_{2};x)=I(G_{1};x)+I(G_{2};x)-1.
\]

\end{proposition}

A finite sequence of real numbers $(a_{0},a_{1},a_{2},...,a_{n})$ is said to
be \textit{unimodal} if there is some $k$, called the \textit{mode} of the
sequence, such that $a_{0}\leq...\leq a_{k-1}\leq a_{k}\geq a_{k+1}\geq...\geq
a_{n}$, and \textit{log-concave} if $a_{i}^{2}\geq a_{i-1}\cdot a_{i+1}$ for
$1\leq i\leq n-1$. It is known that any log-concave sequence of positive
numbers is also unimodal. A polynomial is called \textit{unimodal
(log-concave)} if the sequence of its coefficients is unimodal (log-concave,
respectively). For instance, $I(K_{n}\uplus\left(  \sqcup3K_{7}\right)
;x)=1+(n+21)x+147x^{2}+343x^{3},n\geq1$, is \emph{(a)} log-concave, if
$147^{2}-(n+21)\cdot343\geq0$, i.e., for $1\leq n\leq42$ (e.g., $I(K_{42}%
\uplus\left(  \sqcup3K_{7}\right)  ;x)=1+63x+147x^{2}+343x^{3}$), \emph{(b)}
unimodal, but non-log-concave, whenever $147^{2}-(n+21)\cdot343<0 $ and
$n\leq126$, that is, $43\leq n\leq126$ (for instance, $I(K_{43}\uplus\left(
\sqcup3K_{7}\right)  ;x)=1+64x+147x^{2}+343x^{3}$), \emph{(c)} non-unimodal
for $n\geq127$ (e.g., $I(K_{127}\uplus\left(  \sqcup3K_{7}\right)
;x)=1+148x+147x^{2}+343x^{3}$). The graph $H=(\sqcup3K_{10})\uplus
K_{\underset{120}{\underbrace{3,3,...,3}}}$ is connected and well-covered, but
not very well-covered, and its independence polynomial is unimodal, but not
log-concave: $I(H;x)=1+390x+660x^{2}+\mathbf{1120}x^{3}$. The product of two
polynomials, one log-concave and the other unimodal, is not always
log-concave, for instance, if $G=K_{40}\uplus\left(  \sqcup3K_{7}\right)
,H=K_{110}\uplus\left(  \sqcup3K_{7}\right)  $, then%
\[%
\begin{array}
[c]{l}%
I(G;x)\cdot I(H;x)=\left(  1+61x+147x^{2}+343x^{3}\right)  \left(
1+131x+147x^{2}+343x^{3}\right) \\
=1+192x+8285x^{2}+28910x^{3}+87465x^{4}+100842x^{5}+117649x^{6}.
\end{array}
\]
However, the following result, due to Keilson and Gerber, states that:

\begin{theorem}
\cite{KeilsonGerber}\label{th2} If $P(x)$ is log-concave and $Q(x)$ is
unimodal, then $P(x)\cdot Q(x)$ is unimodal, while the product of two
log-concave polynomials is log-concave.
\end{theorem}

Alavi \textit{et al.} \cite{AlMalSchErdos} showed that for any permutation
$\sigma$ of $\{1,2,...,\alpha\}$ there is a graph $G$ with $\alpha(G)=\alpha$
such that $s_{\sigma(1)}<s_{\sigma(2)}<...<s_{\sigma(\alpha)}$. Nevertheless,
in \cite{AlMalSchErdos} it is stated the following (still open) conjecture:
$I(F;x)$ of any forest $F$ is unimodal.

In \cite{BrownDilNow} it was conjectured that $I(G;x)$ is unimodal for each
well-covered graph $G$. Michael and Traves \cite{MichaelTraves} proved that
this assertion is true for $\alpha(G)\leq3$, but it is false for $4\leq
\alpha(G)\leq7$. In \cite{LevMan6} we showed that for any $\alpha\geq8$, there
exists a connected well-covered graph $G$ with $\alpha(G)=\alpha$, whose
$I(G;x)$ is not unimodal. However, the conjecture of Brown \textit{et al}. is
still open for very well-covered graphs. In \cite{LevMan5} an infinite family
of very well-covered graphs with unimodal independence polynomials is
described. We also showed that $I(G^{\ast};x)$ is unimodal for any $G^{\ast}$
whose skeleton $G$ has $\alpha(G)\leq4$ (see \cite{LevMan5}).

Michael and Traves \cite{MichaelTraves} formulated (and verified for
well-covered graphs with stability numbers $\leq7$) the following so-called
"\textit{roller-coaster}"\ conjecture: for any permutation $\pi$ of the set
$\{\left\lceil \alpha/2\right\rceil ,\left\lceil \alpha/2\right\rceil
+1,...,\alpha\}$, there exists a well-covered graph $G$, with $\alpha
(G)=\alpha$, whose sequence $(s_{0},s_{1},...,s_{\alpha})$ satisfies the
inequalities $s_{\pi\left(  \left\lceil \alpha/2\right\rceil \right)  }%
<s_{\pi\left(  \left\lceil \alpha/2\right\rceil +1\right)  }<...<s_{\pi\left(
\alpha\right)  }$. Recently, Matchett \cite{Matchett} showed that this
conjecture is true for well-covered graphs with stability numbers $\leq11$.

Recall also the following statement, due to Hamidoune.

\begin{theorem}
\cite{Hamidoune}\label{th1} The independence polynomial of a claw-free graph
is log-concave.
\end{theorem}

As a consequence, we deduce that for any $\alpha\geq1$, there exists a tree
$T$, with $\alpha(T)=\alpha$ and whose $I(T;x)$ is log-concave, e.g., the
chordless path $P_{2\alpha}$.

In this paper we show that the independence polynomial of $G^{\ast}$ is
log-concave, whenever: $\alpha(G)\leq3$, or $G^{\ast}$ is a well-covered
spider (i.e., $G=K_{1,n},n\geq1$), or $G^{\ast}$ is a centipede (that is,
$G=P_{n},n\geq1$).

\section{Results}

\begin{lemma}
\label{lem2}If $G$ is a graph of order $n\geq1$ and $\alpha(G)=\alpha$, then
$\alpha\cdot s_{\alpha}\leq n\cdot s_{\alpha-1}$.
\end{lemma}

\begin{proof}
\textbf{\ }Let $H=(\mathcal{A},\mathcal{B},\mathcal{W})$ be the bipartite
graph defined as follows: $X\in\mathcal{A}\Leftrightarrow X$ is a stable set
in $G$ of size $\alpha-1$, then $Y\in\mathcal{B}\Leftrightarrow Y$ is a stable
set in $G$ of size $\alpha(G)$, and $XY\in\mathcal{W}$ $\Leftrightarrow
X\subset Y$ in $G$. Since any $Y\in\mathcal{B}$ has exactly $\alpha(G)$
subsets of size $\alpha-1$, it follows that $\left\vert \mathcal{W}\right\vert
=\alpha\cdot s_{\alpha}$. On the other hand, if $X\in\mathcal{A}$ , then
$\left\vert \{X\cup\{y\}:X\cup\{y\}\in\mathcal{B}\}\right\vert \leq
n-\left\vert X\right\vert =n-\alpha+1$. Hence, any $X\in\mathcal{A}$ has at
most $n-\alpha+1$ neighbors. Consequently, $\left\vert \mathcal{W}\right\vert
=\alpha\cdot s_{\alpha}\leq\left(  n-\alpha+1\right)  \cdot s_{\alpha-1}$, and
this leads to $\alpha\cdot s_{\alpha}\leq n\cdot s_{\alpha-1}$.
\end{proof}

In \cite{LevMan4} it was established the following result:

\begin{theorem}
\label{th4}\cite{LevMan4} If $G$ is a graph of order $n\geq1$ and
$I(G;x)=\sum\limits_{k=0}^{\alpha(G)}s_{k}x^{k}$, then
\[
I(G^{\ast};x)=\sum\limits_{k=0}^{\alpha(G^{\ast})}t_{k}x^{k},\quad t_{k}%
=\sum\limits_{j=0}^{k}s_{j}\cdot\binom{n-j}{n-k},\quad0\leq k\leq
\alpha(G^{\ast})=n.
\]

\end{theorem}

In \cite{LevMan5} it was shown that $I(G^{\ast};x)$ is unimodal for any graph
$G$ with $\alpha(G)\leq4$. Now we partially strengthen this assertion to the
following result.

\begin{theorem}
\label{th5}If $G$ is a graph with $\alpha(G)\leq3$, then $I(G^{*};x)$ is log-concave.
\end{theorem}

\begin{proof}
\textbf{\ }Suppose that\ $\alpha(G)=3$. Then $n=\left\vert V(G)\right\vert
\geq3$ and $I(G;x)=1+nx+s_{2}x^{2}+s_{3}x^{3}$. According to Theorem
\ref{th4}, for $2\leq k\leq n-1$, we obtain: $t_{k}=\binom{n}{k}+n\binom
{n-1}{k-1}+s_{2}\binom{n-2}{k-2}+s_{3}\binom{n-3}{k-3}$. Therefore,
\begin{align*}
t_{k}^{2}-t_{k-1}t_{k+1}  & =A_{0}+n^{2}A_{1}+s_{2}^{2}A_{2}+s_{3}^{2}A_{3}+\\
& nA_{01}+s_{2}A_{02}+s_{3}A_{03}+ns_{2}A_{12}+ns_{3}A_{13}+s_{2}s_{3}A_{23},
\end{align*}
and all $A_{i}\geq0,0\leq i\leq3$, where
\begin{align*}
A_{0}  & =\binom{n}{k}^{2}-\binom{n}{k-1}\binom{n}{k+1},\quad A_{1}%
=\binom{n-1}{k-1}^{2}-\binom{n-1}{k-2}\binom{n-1}{k},\\
A_{2}  & =\binom{n-2}{k-2}^{2}-\binom{n-2}{k-3}\binom{n-2}{k-1},\quad
A_{3}=\binom{n-3}{k-3}^{2}-\binom{n-3}{k-4}\binom{n-3}{k-2},
\end{align*}
because the sequence $\left\{  \binom{n}{k}\right\}  $ is log-concave. Based
on notation $b=\binom{n}{k}^{2}$, we get
\begin{align*}
A_{01}  & =\frac{2k\left(  n+1\right)  b}{n\left(  n-k+1\right)  \left(
k+1\right)  },\quad A_{02}=\frac{2kb\{(k-2)n+2k-1\}}{(k+1)\left(
n-k+1\right)  \left(  n-1\right)  n},\\
A_{03}  & =\frac{2kb(k-1)\{(k-5)n+4k-2\}}{\left(  k+1\right)  \left(
n-k+1\right)  \left(  n-2\right)  \left(  n-1\right)  n},\quad A_{12}%
=\frac{2kb(k-1)}{n(n-1)\left(  n-k+1\right)  },\\
A_{13}  & =\frac{2kb\left(  k-1\right)  \{(k-3)n+k\}}{\left(  n-2\right)
\left(  n-1\right)  n^{2}(n-k+1)},\quad A_{23}=\frac{2k^{2}b\left(
k-1\right)  \left(  k-2\right)  }{\left(  n-2\right)  \left(  n-1\right)
n^{2}\left(  n-k+1\right)  },
\end{align*}
and all $A_{ij}\geq0$ for $k\geq5$. Hence, we must check that $t_{k}%
^{2}-t_{k-1}t_{k+1}\geq0$ for $k\in\{1,2,3,4\}$. By Theorem \ref{th4}, we
obtain:
\[
t_{0}=1,t_{1}=2n,t_{2}=3n(n-1)/2+s_{2},t_{3}=\allowbreak\frac{2n(n-1)(n-2)}%
{3}+(n-2)s_{2}+s_{3},
\]%
\[
t_{4}=\frac{5}{24}\left(  n-3\right)  n\left(  n-1\right)  \left(  n-2\right)
+\frac{1}{2}s_{2}\left(  n-2\right)  \left(  n-3\right)  +s_{3}n-3s_{3},
\]%
\[
t_{5}=\left(  n-4\right)  \left(  n-3\right)  \left[  \frac{1}{20}n\left(
n-1\right)  \left(  n-2\right)  +\frac{1}{6}s_{2}\left(  n-2\right)  +\frac
{1}{2}s_{3}\right]  .
\]

Consequently, it follows $t_{1}^{2}-t_{0}t_{1}=\left(  n^{2}+2\left(
2n^{2}-s_{2}\right)  \right)  /2>0$. We also deduce
\[
t_{2}^{2}-t_{1}t_{3}=\frac{1}{12}\left(  11n+5\right)  (n-1)n^{2}+s_{2}%
^{2}+ns_{2}+n\left(  ns_{2}-2s_{3}\right)  \geq0,
\]
since $3s_{3}\leq ns_{2}$ is true according to Lemma \ref{lem2}.

Now, simple calculations lead us to%
\[%
\begin{array}
[c]{l}%
144\left(  t_{3}^{2}-t_{2}t_{4}\right)  =\left(  19n+7\right)  n^{2}\left(
n-1\right)  ^{2}\left(  n-2\right)  +\left(  54n+30\right)  n\left(
n-1\right)  \left(  n-2\right)  s_{2}\\
-24\left(  n-11\right)  n\left(  n-1\right)  s_{3}+72n\left(  n-3\right)
s_{2}^{2}+144(s_{3}^{2}+(n-1)s_{2}s_{3}+s_{2}^{2}).
\end{array}
\]
Let us notice that $n\left(  n-1\right)  \left(  \left(  54n+30\right)
\left(  n-2\right)  s_{2}-24\left(  n-11\right)  s_{3}\right)  \geq0$, because
Lemma \ref{lem2} implies the inequality $54ns_{2}\geq24s_{3}$. Hence, we infer
that $t_{3}^{2}-t_{2}t_{4}\geq0$, whenever $n\geq3$.

Further, we have%
\[%
\begin{array}
[c]{l}%
2880\left(  t_{4}^{2}-t_{3}t_{5}\right)  =\left(  29n^{8}-252n^{7}%
-108n^{2}+818n^{6}+12n^{3}-1200n^{5}+701n^{4}\right)  +\\
+\left(  672n+2680n^{4}-2520n^{3}+64n^{2}+136n^{6}-1032n^{5}\right)  s_{2}+\\
+\left(  -3840n^{3}+8400n^{2}-4896n+96n^{5}+240n^{4}\right)  s_{3}+\\
+\left(  240n^{4}-1920n^{3}+5520n^{2}-6720n+2880\right)  s_{2}^{2}+\\
+\left(  10\,560n-5760n^{2}+960n^{3}-5760\right)  s_{3}s_{2}+\left(
8640+1440n^{2}-7200n\right)  s_{3}^{2}\\
=\left(  29n+9\right)  n^{2}\left(  n-1\right)  ^{2}\left(  n-2\right)
^{2}\left(  n-3\right)  +\\
+\left(  136n+56\right)  n\left(  n-1\right)  \left(  n-2\right)  ^{2}\left(
n-3\right)  s_{2}+\\
+\left(  96n+816\right)  n\left(  n-1\right)  \left(  n-2\right)  \left(
n-3\right)  s_{3}+240\left(  n-1\right)  \left(  n-2\right)  ^{2}\left(
n-3\right)  s_{2}^{2}+\\
+960\left(  n-1\right)  \left(  n-2\right)  \left(  n-3\right)  s_{2}%
s_{3}+1440\left(  n-2\right)  \left(  n-3\right)  s_{3}^{2}\geq0.
\end{array}
\]
Consequently, $t_{k}^{2}-t_{k-1}t_{k+1}\geq0$, for $1\leq k\leq n-1$, i.e.,
$I(G^{\ast};x)$ is log-concave.

The log-concavity for the cases $\alpha(G)\in\{1,2\}$ can be validated in a
similar way, by observing that either $s_{2}=s_{3}=0$ or only $s_{3}=0$.
\end{proof}

Since $\alpha(K_{1,n})=n,\alpha(P_{n})=\left\lceil n/2\right\rceil $, Theorem
\ref{th5} is not useful in proving that $I(K_{1,n}^{\ast};x)$, $I(P_{n}^{\ast
};x)$ are log-concave, as soon as $n$ is sufficiently large. In \cite{LevMan2}%
, \cite{LevMan3} we proved that $I(K_{1,n}^{\ast};x)$, $I(W_{n};x)$ are
unimodal. Here we are strengthening these results.

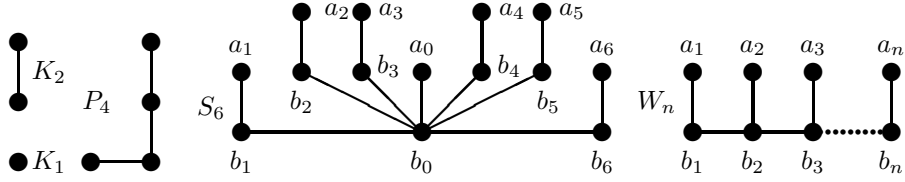
\begin{figure}[h]
\setlength{\unitlength}{0.8cm} \begin{picture}(5,2.5)\thicklines
\put(-4.7,0){\circle*{0.29}}
\put(-4.2,0){\makebox(0,0){$K_{1}$}}
\put(-4.7,1){\circle*{0.29}}
\put(-4.7,2){\circle*{0.29}}
\put(-4.7,1){\line(0,1){1}}
\put(-4.2,1.5){\makebox(0,0){$K_{2}$}}
\multiput(-3.5,0)(1,0){2}{\circle*{0.29}}
\multiput(-2.5,1)(0,1){2}{\circle*{0.29}}
\put(-3.5,0){\line(1,0){1}}
\multiput(-2.5,0)(0,1){2}{\line(0,1){1}}
\put(-3.4,1){\makebox(0,0){$P_{4}$}}
\multiput(-1,0.5)(3,0){3}{\circle*{0.29}}
\multiput(-1,1.5)(1,0){7}{\circle*{0.29}}
\multiput(0,2.5)(1,0){2}{\circle*{0.29}}
\multiput(3,2.5)(1,0){2}{\circle*{0.29}}
\multiput(-1,0.5)(1,0){6}{\line(1,0){1}}
\multiput(-1,0.5)(3,0){3}{\line(0,1){1}}
\multiput(0,1.5)(1,0){2}{\line(0,1){1}}
\multiput(3,1.5)(1,0){2}{\line(0,1){1}}
\put(2,0.5){\line(-2,1){2}}
\put(2,0.5){\line(-1,1){1}}
\put(2,0.5){\line(1,1){1}}
\put(2,0.5){\line(2,1){2}}
\put(2,0){\makebox(0,0){$b_{0}$}}
\put(2,1.9){\makebox(0,0){$a_{0}$}}
\put(-1,0){\makebox(0,0){$b_{1}$}}
\put(-1,1.9){\makebox(0,0){$a_{1}$}}
\put(0,1){\makebox(0,0){$b_{2}$}}
\put(0.6,2.5){\makebox(0,0){$a_{2}$}}
\put(1.45,1.6){\makebox(0,0){$b_{3}$}}
\put(1.5,2.5){\makebox(0,0){$a_{3}$}}
\put(3.45,1.6){\makebox(0,0){$b_{4}$}}
\put(3.5,2.5){\makebox(0,0){$a_{4}$}}
\put(4.1,1){\makebox(0,0){$b_{5}$}}
\put(4.5,2.5){\makebox(0,0){$a_{5}$}}
\put(5,0){\makebox(0,0){$b_{6}$}}
\put(5,1.9){\makebox(0,0){$a_{6}$}}
\put(-1.5,0.9){\makebox(0,0){$S_{6}$}}
\multiput(6.5,0.5)(1,0){3}{\circle*{0.29}}
\multiput(6.5,1.5)(1,0){3}{\circle*{0.29}}
\put(9.8,0.5){\circle*{0.29}}
\put(9.8,1.5){\circle*{0.29}}
\put(6.5,0.5){\line(0,1){1}}
\put(7.5,0.5){\line(0,1){1}}
\put(8.5,0.5){\line(0,1){1}}
\put(6.5,0.5){\line(1,0){1}}
\put(7.5,0.5){\line(1,0){1}}
\put(9.8,0.5){\line(0,1){1}}
\multiput(8.5,0.5)(0.125,0){10}{\circle*{0.07}}
\put(6.5,0){\makebox(0,0){$b_{1}$}}
\put(6.5,1.9){\makebox(0,0){$a_{1}$}}
\put(7.5,0){\makebox(0,0){$b_{2}$}}
\put(7.5,1.9){\makebox(0,0){$a_{2}$}}
\put(8.5,0){\makebox(0,0){$b_{3}$}}
\put(8.5,1.9){\makebox(0,0){$a_{3}$}}
\put(9.8,0){\makebox(0,0){$b_{n}$}}
\put(9.8,1.9){\makebox(0,0){$a_{n}$}}
\put(5.9,1){\makebox(0,0){$W_{n}$}}
\end{picture}
\caption{Well-covered spiders:\ $K_{1},K_{2},P_{4},S_{6}$, and the centipede
$W_{n}$.}%
\label{fig999}%
\end{figure}

The well-covered spider $S_{n},n\geq2$, has $n$ vertices of degree $2$, one
vertex of degree $n+1$, and $n+1$ vertices of degree $1$ (see Figure
\ref{fig999}). In fact, it is easy to see that $S_{n}=K_{1,n}^{\ast},n\geq2$.

\begin{proposition}
\label{prop1}\cite{LevMan3} The independence polynomial of any well-covered
spider is unimodal, moreover, $I(S_{n};x)=(1+x)\cdot\sum\limits_{k=0}%
^{n}\left[  \binom{n}{k}\cdot2^{k}+\binom{n-1}{k-1}\right]  \cdot x^{k}%
,n\geq2$, and its mode is unique and equals $1+\left(  n-1\right)
\operatorname{mod}3+2\left(  \left\lceil n/3\right\rceil -1\right)  $.
\end{proposition}

In \cite{BrownDilNow} it was shown that $I(G;x)$ of any graph $G$ with
$\alpha(G)=2$ has real roots, and, hence, it is log-concave, according to
Newton's theorem (stating that if a polynomial with positive coefficients has
only real roots, then its coefficients form a log-concave sequence). However,
Newton's theorem is not useful in solving the conjecture of Alavi \textit{et
al}., even for the particular case of very well-covered\emph{\ }trees, since,
for instance, $I(S_{3};x)=1+8x+21x^{2}+23x^{3}+9x^{4}$ has non-real roots.

\begin{theorem}
\label{th3}The independence polynomial of any well-covered spider is log-concave.
\end{theorem}

\begin{proof}
Since $I(G;x)$ is log-concave for any graph $G$ with $\alpha(G)\leq2$, we
consider only well-covered spiders $S_{n}$ with $n\geq2$. According to
Proposition \ref{prop1},
\[
I(S_{n};x)=(1+x)\cdot\sum\limits_{k=0}^{n}\left[  \binom{n}{k}\cdot
2^{k}+\binom{n-1}{k-1}\right]  \cdot x^{k}=(1+x)\cdot P(x).
\]
It is sufficient to prove that $P(x)$ is log-concave, because, further,
Theorem \ref{th2} implies that $I(S_{n};x)$ is log-concave, as well. Let us
denote $c_{k}=\binom{n}{k}\cdot2^{k}+\binom{n-1}{k-1},0\leq k\leq n$.

Firstly, we notice that $c_{1}^{2}-c_{0}\cdot c_{2}=(2n+1)(n+2)>0$. Further,
for $2\leq k\leq n-1$, we obtain that:
\begin{align*}
c_{k}^{2}-c_{k-1}\cdot c_{k+1}  & =\left[  \binom{n-1}{k-1}^{2}-\binom
{n-1}{k-2}\binom{n-1}{k}\right]  +\\
& +\binom{n}{k}^{2}\frac{n(2n+2)2^{k}-k^{2}\left(  n+3\right)  +k(k^{2}%
+7n+4)}{n\left(  k+1\right)  \left(  n-k+1\right)  \cdot2^{1-k}}.
\end{align*}
Clearly, $\binom{n-1}{k-1}^{2}-\binom{n-1}{k-2}\binom{n-1}{k}\geq0$, since the
sequence of binomial coefficients is log-concave, and $n(2n+2)2^{k}%
-k^{2}\left(  n+3\right)  \geq0$, because $n\cdot2^{k}\geq k^{2}$ holds for
any $k\in\{2,....,n-1\}$. Thus, $c_{k}^{2}-c_{k-1}\cdot c_{k+1}\geq0$, for any
$k\in\{1,2,...,n-1\}$.
\end{proof}

The \textit{edge-join} of two disjoint graphs $G_{1},G_{2}$, is the graph
$G_{1}\circleddash G_{2}$ obtained by adding an edge joining a vertex from
$G_{1}$ to a vertex from $G_{2}$. If both vertices are of degree at least two,
then $G_{1}\circleddash G_{2}$ is an \textit{internal edge-join} of
$G_{1},G_{2}$. By $\bigtriangleup_{n}$ we mean the graph $\circleddash
nK_{3}=(\circleddash(n-1)K_{3})\circleddash K_{3},n\geq1$ (see Figure
\ref{fig22}). \begin{figure}[h]
\setlength{\unitlength}{1cm}\begin{picture}(5,1)\thicklines
\multiput(-2,0)(1,0){6}{\circle*{0.29}}
\multiput(-1.5,1)(2,0){3}{\circle*{0.29}}
\put(5,0){\circle*{0.29}}
\put(6,0){\circle*{0.29}}
\put(5.5,1){\circle*{0.29}}
\put(-2,0){\line(1,0){5}}
\put(-2,0){\line(1,2){0.5}}
\put(0,0){\line(1,2){0.5}}
\put(2,0){\line(1,2){0.5}}
\put(5,0){\line(1,2){0.5}}
\put(-1,0){\line(-1,2){0.5}}
\put(1,0){\line(-1,2){0.5}}
\put(3,0){\line(-1,2){0.5}}
\put(6,0){\line(-1,2){0.5}}
\put(5,0){\line(1,0){1}}
\multiput(3,0)(0.125,0){16}{\circle*{0.07}}
\put(-2.4,0){\makebox(0,0){$v_{1}$}}
\put(-0.9,0.3){\makebox(0,0){$v_{2}$}}
\put(-0.1,0.3){\makebox(0,0){$v_{4}$}}
\put(1.1,0.3){\makebox(0,0){$v_{5}$}}
\put(1.9,0.3){\makebox(0,0){$v_{7}$}}
\put(3.1,0.3){\makebox(0,0){$v_{8}$}}
\put(4.5,0.3){\makebox(0,0){$v_{3n-2}$}}
\put(6.65,0){\makebox(0,0){$v_{3n-1}$}}
\put(-1.85,1){\makebox(0,0){$v_{3}$}}
\put(0.85,1){\makebox(0,0){$v_{6}$}}
\put(2.85,1){\makebox(0,0){$v_{9}$}}
\put(6,1){\makebox(0,0){$v_{3n}$}}
\end{picture}
\caption{The graph $\bigtriangleup_{n}=(\circleddash(n-1)K_{3})\circleddash
K_{3}$.}%
\label{fig22}%
\end{figure}

In \cite{FinHarNow} it is shown that apart from $K_{1}$ and $C_{7}$, any
connected well-covered graph of girth $\geq6$ equals $G^{\ast}$ for some graph
$G$, e.g., every well-covered tree equals $T^{\ast}$ for some tree $T$ (see
also \cite{Ravindra}). Thus, a tree $T\neq K_{1}$ could be only very well-covered.

\begin{theorem}
\cite{LevMan1} A tree $T$ is well-covered if and only if $T$ is a well-covered
spider, or $T$ is the internal edge-join of a number of well-covered spiders.
\end{theorem}

A \textit{centipede }is a well-covered tree\textit{\ }defined by $W_{n}%
=P_{n}^{\ast},n\geq1$ (see\textit{\ }Figure \ref{fig999}). For example,
$W_{1}=K_{2},W_{2}=P_{4},W_{3}=S_{2}$.

\begin{theorem}
\label{th6}The independence polynomial of any centipede is log-concave.
\end{theorem}

\begin{proof}
\textbf{\ }We show, by induction on $n\geq1$, that
\[
I(W_{2n};x)=(1+x)^{n}\cdot I(\bigtriangleup_{n};x),\quad I(W_{2n+1}%
;x)=(1+x)^{n}\cdot I(\bigtriangleup_{n}\circleddash K_{2};x),
\]
(for another proof of these equalities, see \cite{LevMan3}).

For $n=1$, the assertion is true, because
\begin{align*}
I(W_{2};x)  & =1+4x+3x^{2}=(1+x)(1+3x)=(1+x)\cdot I(\bigtriangleup_{1};x),\\
I(W_{3};x)  & =1+6x+10x^{2}+5x^{3}=(1+x)\cdot I(\bigtriangleup_{1}\circleddash
K_{2};x).
\end{align*}

Assume that the formulae are true for $k\leq2n+1$. By Proposition \ref{prop4},
we get:
\begin{align*}
I(W_{2n+2};x)  & =I(W_{2n+2}-b_{2n+1};x)+x\cdot I(W_{2n+2}-N[b_{2n+1}];x)\\
& =\left(  1+x\right)  \left(  1+2x\right)  \cdot I(W_{2n};x)+x\left(
1+x\right)  ^{2}\cdot I(W_{2n-1};x)\\
& =\left(  1+x\right)  ^{n+1}\cdot\left\{  I(K_{2};x)\cdot I(\bigtriangleup
_{n};x)+x\cdot I(\bigtriangleup_{n-1}\circleddash K_{2};x))\right\}  .
\end{align*}
On the other hand, if $v$ is the vertex of degree $3$ in the last triangle of
$\bigtriangleup_{n+1}$(see Figure \ref{fig1234}(\textit{a})), then
$I(\bigtriangleup_{n+1};x)=I(K_{2};x)I(\bigtriangleup_{n};x)+xI(\bigtriangleup
_{n-1}\circleddash K_{2};x))$, according to Proposition \ref{prop4}.
\begin{figure}[h]
\setlength{\unitlength}{1cm}\begin{picture}(5,2.5)\thicklines
\multiput(-3,1.5)(1,0){4}{\circle*{0.29}}
\multiput(-2.5,2.5)(2,0){2}{\circle*{0.29}}
\put(3,1.5){\circle*{0.29}}
\put(2,1.5){\circle*{0.29}}
\put(2.5,2.5){\circle*{0.29}}
\put(-3,1.5){\line(1,0){3}}
\put(-3,1.5){\line(1,2){0.5}}
\put(-1,1.5){\line(1,2){0.5}}
\put(2,1.5){\line(1,2){0.5}}
\put(-2,1.5){\line(-1,2){0.5}}
\put(0,1.5){\line(-1,2){0.5}}
\put(3,1.5){\line(-1,2){0.5}}
\put(2,1.5){\line(1,0){1}}
\multiput(0,1.5)(0.125,0){16}{\circle*{0.07}}
\put(1.8,1.85){\makebox(0,0){$v$}}
\put(-3.5,2){\makebox(0,0){$(a)$}}
\multiput(4.3,1.5)(1,0){2}{\circle*{0.29}}
\multiput(4.3,2.5)(1,0){2}{\circle*{0.29}}
\multiput(6.3,1.5)(1.2,0){2}{\circle*{0.29}}
\multiput(6.3,2.5)(1.2,0){2}{\circle*{0.29}}
\put(4.3,1.5){\line(1,0){1}}
\put(4.3,1.5){\line(0,1){1}}
\put(5.3,1.5){\line(0,1){1}}
\put(6.3,1.5){\line(1,0){1.2}}
\put(6.3,1.5){\line(0,1){1}}
\put(7.5,1.5){\line(0,1){1}}
\multiput(5.3,1.5)(0.125,0){9}{\circle*{0.07}}
\put(4,2.5){\makebox(0,0){$a_{1}$}}
\put(4,1.5){\makebox(0,0){$b_{1}$}}
\put(6.9,2.3){\makebox(0,0){$a_{2n+1}$}}
\put(6.9,1.8){\makebox(0,0){$b_{2n+1}$}}
\put(8.15,2.5){\makebox(0,0){$a_{2n+2}$}}
\put(8.15,1.5){\makebox(0,0){$b_{2n+2}$}}
\multiput(-3,0)(1,0){4}{\circle*{0.29}}
\multiput(-2.5,1)(2,0){2}{\circle*{0.29}}
\put(1,0){\circle*{0.29}}
\put(2,0){\circle*{0.29}}
\put(3,0){\circle*{0.29}}
\put(3,1){\circle*{0.29}}
\put(1.5,1){\circle*{0.29}}
\put(-3,0){\line(1,0){3}}
\put(-3,0){\line(1,2){0.5}}
\put(-1,0){\line(1,2){0.5}}
\put(1,0){\line(1,2){0.5}}
\put(-2,0){\line(-1,2){0.5}}
\put(0,0){\line(-1,2){0.5}}
\put(2,0){\line(-1,2){0.5}}
\put(1,0){\line(1,0){2}}
\put(3,0){\line(0,1){1}}
\multiput(0,0)(0.125,0){9}{\circle*{0.07}}
\put(2,0.3){\makebox(0,0){$v$}}
\put(-3.5,0.5){\makebox(0,0){$(b)$}}
\multiput(4.3,0)(1,0){2}{\circle*{0.29}}
\multiput(4.3,1)(1,0){2}{\circle*{0.29}}
\multiput(6.3,0)(1.2,0){2}{\circle*{0.29}}
\multiput(6.3,1)(1.2,0){2}{\circle*{0.29}}
\put(4.3,0){\line(1,0){1}}
\put(4.3,0){\line(0,1){1}}
\put(5.3,0){\line(0,1){1}}
\put(6.3,0){\line(1,0){1.2}}
\put(6.3,0){\line(0,1){1}}
\put(7.5,0){\line(0,1){1}}
\multiput(5.3,0)(0.125,0){9}{\circle*{0.07}}
\put(4,1){\makebox(0,0){$a_{1}$}}
\put(4,0){\makebox(0,0){$b_{1}$}}
\put(6.9,0.8){\makebox(0,0){$a_{2n+2}$}}
\put(6.9,0.3){\makebox(0,0){$b_{2n+2}$}}
\put(8.15,1){\makebox(0,0){$a_{2n+3}$}}
\put(8.15,0){\makebox(0,0){$b_{2n+3}$}}
\end{picture}
\caption{The graphs: (\textit{a}) $\bigtriangleup_{n+1}$ and $W_{2n+2}$;
(\textit{b}) $\bigtriangleup_{n+1}\circleddash K_{2}$ and $W_{2n+3}$.}%
\label{fig1234}%
\end{figure}
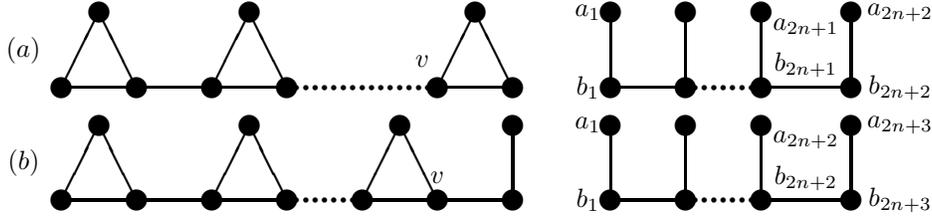In other words, $I(W_{2n+2};x)=\left(  1+x\right)  ^{n+1}\cdot
I(\bigtriangleup_{n+1};x)$.

Similarly, again by Proposition \ref{prop4}, we obtain:
\begin{align*}
I(W_{2n+3};x)  & =I(W_{2n+3}-b_{2n+2};x)+x\cdot I(W_{2n+3}-N[b_{2n+2}];x)\\
& =\left(  1+x\right)  \left(  1+2x\right)  \cdot I(W_{2n+1};x)+x\left(
1+x\right)  ^{2}\cdot I(W_{2n};x)\\
& =\left(  1+x\right)  ^{n+1}\left\{  I(K_{2};x)\cdot I(\bigtriangleup
_{n}\circleddash K_{2};x)+x\left(  1+x\right)  \cdot I(\bigtriangleup
_{n};x))\right\}  .
\end{align*}

On the other hand, if $v$ is the vertex of degree $3$ belonging to the last
triangle of $\bigtriangleup_{n+1}\circleddash K_{2}$ (see Figure
\ref{fig1234}(\textit{b})) and adjacent to one of the vertices of $K_{2}$, we
have
\begin{align*}
I(\bigtriangleup_{n+1}\circleddash K_{2};x)  & =I(\bigtriangleup
_{n+1}\circleddash K_{2}-v;x)+xI(\bigtriangleup_{n+1}\circleddash
K_{2}-N[v];x)\\
& =I(K_{2};x)\cdot I(\bigtriangleup_{n}\circleddash K_{2};x)+x\left(
1+x\right)  \cdot I(\bigtriangleup_{n};x)).
\end{align*}
In other words,%
\[
I(W_{2n+3};x)=\left(  1+x\right)  ^{n+1}\cdot I(\bigtriangleup_{n+1}%
\circleddash K_{2};x).
\]

While Theorem \ref{th1} assures that $I(\bigtriangleup_{n};x),I(\bigtriangleup
_{n}\circleddash K_{2};x)$ are log-concave, finally Theorem \ref{th2} implies
that $I(W_{n};x)$ is log-concave, as claimed.
\end{proof}

\begin{corollary}
\emph{(i)} If the graph $H$ has as connected components well-covered

spiders/centipedes and/or graphs with stability number $\leq2$, and/or

claw-free graphs, and/or graphs that may be represented as $G^{\ast}$ whose

$G$ has $\alpha(G)\leq3$, then its independence polynomial $I(H;x)$ is log-concave.

\emph{(ii)} If $H_{n}\in\{S_{n},W_{n}\}$, then the independence polynomial of
$\uplus mH_{n}$ is

log-concave, for any $m\geq2,n\geq1$.
\end{corollary}

\begin{proof}
\textbf{\ }\emph{(i)} Let $G_{i},1\leq i\leq m$, be the connected components
of $G$. According to Theorems \ref{th3}, \ref{th6}, \ref{th5} and \ref{th1},
any $I(G_{i};x)$ is log-concave. Further, Theorem \ref{th2} implies that
$I(G;x)$ is also log-concave, as $I(G;x)=I(G_{1};x)\cdot...\cdot I(G_{m};x)$.

\emph{(ii)} Since $I(H_{n};x)$ is log-concave, and $I(\uplus mH_{n};x)=m\cdot
I(H_{n};x)-(m-1)$, it follows that $I(\uplus mH_{n};x)$ is log-concave, as well.
\end{proof}

\section{Conclusions}

In this paper we showed that for any $\alpha$, there is a very well-covered
tree $T$ with $\alpha(T)=\alpha$, whose independence polynomial $I(T;x)$ is
log-concave. We conjecture that the independence polynomial of any\emph{\ }%
(well-covered)\emph{\ }forest is log-concave.

\begin{figure}[h]
\setlength{\unitlength}{1cm}\begin{picture}(5,2.2)\thicklines
\multiput(-1,0)(1,0){4}{\circle*{0.29}}
\multiput(-1,1)(1,0){4}{\circle*{0.29}}
\multiput(-1,2)(1,0){2}{\circle*{0.29}}
\put(-1,0){\line(1,0){1}}
\put(-1,1){\line(1,0){3}}
\put(-1,2){\line(1,0){1}}
\put(0,0){\line(0,1){2}}
\put(1,0){\line(0,1){1}}
\put(2,0){\line(0,1){1}}
\put(-1.7,1){\makebox(0,0){$T_{1}$}}
\multiput(4,0)(1,0){3}{\circle*{0.29}}
\multiput(4,1)(1,0){3}{\circle*{0.29}}
\multiput(4,2)(1,0){2}{\circle*{0.29}}
\put(4,0){\line(1,0){1}}
\put(4,1){\line(1,0){2}}
\put(4,2){\line(1,0){1}}
\put(5,0){\line(0,1){2}}
\put(6,0){\line(0,1){1}}
\put(3.3,1){\makebox(0,0){$T_{2}$}}
\end{picture}
\caption{Two (very) well-covered trees.}%
\label{fig5}%
\end{figure}
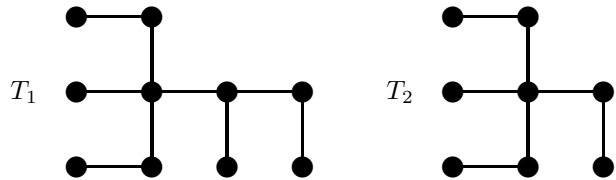

In 1990, Hamidoune \cite{Hamidoune} conjectured that the independence
polynomial of any claw-free graph has only real roots. Recently, Chudnovsky
and Seymour \cite{ChudSeymour} validated this conjecture. Consequently,
$I(P_{n};x)$ has all the roots real. Moreover, the roots of $I(W_{n};x)$ are
real (see the proof\emph{\ }of Theorem \ref{th6}).

For general (very well-covered) spiders/trees the structure of the roots of
the independence polynomial is more complicated. For instance, the
independence polynomial of the claw graph $I(K_{1,3};x)=1+4x+3x^{2}+x^{3}$ has
non-real roots. Figure \ref{fig5} provides us with some more examples:
\begin{align*}
I(T_{1};x)  & =(1+x)^{2}(1+2x)(1+6x+7x^{2}),\\
I(T_{2};x)  & =(1+x)(1+7x+14x^{2}+9x^{3}),
\end{align*}
where only $I(T_{1};x)$ has all the\emph{\ }roots real. It seems to be
interesting to characterize (well-covered) trees whose independence
polynomials have only real roots.

\end{document}